\documentclass{amsart}
\usepackage{amscd} 
\usepackage{amsfonts} 
\usepackage{amssymb} 
\usepackage{latexsym} 
\input{diagrams}
\newcommand{\ncm}{\newcommand}

\def\N{\mathbb{N}\,}

\def\C{\mathbb{C}\,}

\def\MC{\mathcal{C}}

\newtheorem{theorem}{Theorem}[section]
\newtheorem{prop}[theorem]{Proposition}
\newtheorem{lemma}[theorem]{Lemma}

\newtheorem{lem&def}[theorem]{Lemma \& Definition}

\ncm{\End}{\mbox{\rm End}\,}

\def\Hom{\mbox{\rm Hom}\,}

\def\Im{\mbox{\rm Im}\,}

\def\id{\mbox{\rm id}}

\def\into{\hookrightarrow}
\def\to{\rightarrow}

\def\o{\otimes}    %tensor product 
\def\bra{\langle}
\def\ket{\rangle}

\ncm{\rarr}[1]{\stackrel{#1}{\longrightarrow}}
\ncm{\larr}[1]{\stackrel{#1}{\longleftarrow}}
\def\cop{\Delta}

\def\eps{\varepsilon}

\def\du1{\hat 1}

\def\-1{_{(-1)}}
\def\0{_{(0)}}
\def\1{_{(1)}}
\def\2{_{(2)}}
\def\3{_{(3)}}

\def\|{\, | \,}

\def\du1{\hat 1}
\def\lact{\triangleright}

\begin{document}

\title{Depth Two and the Galois Coring}
\author{Lars Kadison}
\address{Matematiska Institutionen \\ G{\" o}teborg 
University \\ 
S-412 96 G{\" o}teborg, Sweden} 
\email{lkadison@c2i.net} 
\date{}
\thanks{The author thanks   P.\ Hajac and K.~Szlachanyi 
for encouraging discussions.}
\date{} 

\begin{abstract} 
We study the cyclic module ${}_SR$ for a ring extension $A \| B$ with centralizer $R$
and bimodule endomorphism ring $S = \End {}_BA_B$.  We show that if
$A \| B$ is an H-separable Hopf subalgebra, then $B$ is a normal Hopf subalgebra
of $A$.  We observe from \cite{KN} and \cite{KS} depth two in the role of noncommutative
normality (as in field theory) in a depth two separable Frobenius characterization of semisimple-Hopf-Galois
extensions.  We prove that a depth two extension has a Galois $A$-coring structure on $A \o_R T$
where $T$ is the right $R$-bialgebroid dual to $S$.  
\end{abstract} 
\maketitle

\section{Introduction}

Depth two extensions have their origins in finite depth $II_1$ subfactors.
An inclusion of finite-dimensional C$^*$-algebras $B \subseteq A$ can be recorded as a bicolored weighted
multigraph called a Bratteli inclusion diagram:  
the number of edges between a black dot representing an (isomorphism class of a)
simple module $V$ of $A$
and a white dot representing a simple module $W$ of $B$ is $\dim {\rm Hom}_B(V,W)$.  This can be recorded in an inclusion matrix
of non-negative integers, which corresponds to an induction-restriction table of
irreducible characters of a  subgroup pair $H < G$ if $A = \C G$ and $B=\C H$ are the group algebras.  

If we define the basic construction of a semisimple $\C$-algebra pair $B \subseteq A$ to be the endomorphism algebra
$\End A_B$ of intertwiners, we note that $B$ and $\mathcal{E} := \End A_B$ are Morita equivalent via bimodule ${}_{\mathcal{E}}A_B$ 
whence the inclusion diagram of the left multiplication inclusion $A \into \End A_B$ is reflection of the 
diagram of $B \subseteq A$.  Beginning with a subfactor $N \subseteq M$, we build the Jones tower using the basic construction
$$ N \subseteq M \subseteq M_1 \subseteq M_2 \subseteq \cdots $$
where $M_{i+1} = \End {M_i}_{M_{i-1}}$, then the derived tower of centralizers or relative commutants
are f.d.\ C$^*$-algebras,
$$ C_N(N) \subseteq C_M(N) \subseteq C_{M_1}(N) \subseteq C_{M_2}(N) \subseteq \cdots $$
The subfactor $N \subseteq M$ has finite depth if the inclusion diagrams of the derived tower stop growing and begin
reflecting at some point, depth $n$ where counting begins with $0$. For example, the subfactor has depth two
if $C_{M_2}(N)$ is isomorphic to the basic construction of $C_M(N) \subseteq C_{M_1}(N) $.

Split depth two extensions are automatically finitely generated (f.g.) projective
(see Prop.~\ref{prop-fgp} below) - a depth two subfactor has finite index and is a Frobenius extension
via an ``orthonormal'' basis w.r.t.\ its conditional expectation.  
In \cite{KN} we define a depth two Frobenius extension $N \subseteq M$ to
be one whose dual bases may be chosen from $C_{M_1}(N)$, and provide a depth two characterization of
certain irreducible semisimple-Hopf-Galois extensions: see Theorem~\ref{th-char} below where depth two
fills the role of normal extension in the theorem, separable + normal = Galois.  In \cite{KS} we widen
the definition of depth two to ring extensions, and show that depth two extensions have (generalized Lu) bialgebroids attached
to them that are dual w.r.t.\ the centralizer and act on the overalgebra and an endomorphism ring. 
In Theorem~\ref{th-Galois coring}
below we take the coaction point of view and show that a depth two extension
forms a Galois coring in a natural way.

The basic set-up throughout this paper is the following. Let
$A \| B$ be a ring extension with centralizer denoted by $R := C_A(B) = A^B$,
bimodule endomorphism ring $S := \End {}_BA_B$ and $B$-central
tensor-square $T := (A \o_B A)^B$. 
  $T$ has a ring structure induced
from $T \cong \End {}_A A \o_B A_A$ given by 
$$tt' = {t'}^1 t^1 \o t^2 {t'}^2 \ \ \ 1_T = 1 \o 1,$$
where $t = t^1 \o t^2 \in T$ uses a Sweedler notation and suppresses
a possible summation over simple tensors.  Let $\lambda: A \into \End A_B$
denote left multiplication and $\rho: A \into \End {}_BA$ denote right multiplication.
Let $\mathcal{E}$ denote $\End A_B$ and
note that $S \subseteq \mathcal{E}$, where $\mathcal{E}_S$ will denote
 the
natural module.  Note that $\lambda$ restricts to $R \into S$ and $\rho$ restricts
to $R \into S^{\rm op}$. Let $\mathcal{Z}(B)$ denote the center of any ring $B$.

We introduce a handy notion in our field, the notion of an arbitrary bimodule
being centrally projective with respect to a canonical bimodule.  We say that
a bimodule ${}_AM_B$, where $A$ and $B$ are two arbitrary rings,
is \textit{centrally projective w.r.t.\ a bimodule} ${}_AN_B$, if ${}_AM_B$ is
isomorphic to a direct summand of a finite direct sum of $N$; in symbols,
if ${}_AM_B \oplus * \cong \oplus^n {}_AN_B$.  This covers the usual notion
of centrally projective $A$-$A$-bimodule $P$ where the canonical $A$-$A$-bimodule
is understood by default to be the natural bimodule $A$ itself. 
 
In \cite{KK} we studied the cyclic module $R_T$ given by a generalized
Miyashita-Ulbrich action: $r \cdot t = t^1 r t^2$ where $r \in R$
and $t \in T$ and a ternary product mapping $\gamma: R \o_T (A \o_B A) \to
A$ given by $r \o a \o a' \mapsto ara'$.  The 
multiplication mapping $\mu: A \o_B A \to A$,
$a \o a \mapsto aa'$ factors through  
 a canonical epi given by $a \o a' \mapsto 1 \o a \o a'$ and the ternary product mapping. 
We showed in \cite{KK} that  $\gamma$ is an isomorphism for depth two
or separable extensions, $\End R_T \cong \mathcal{Z}(A)$ if D2,
and characterizes separable or H-separable extensions if additionally
$R_T$ is f.g.\ projective or a generator module, respectively.  

In this paper we will study the cyclic module
${}_SR$, also an interesting module.  It is just given by evaluation:
\begin{equation}
\label{eq: S-action}
 \alpha \cdot r = \alpha(r) \in R \ \ \ \ (\alpha \in S,\ r \in R )
\end{equation}
which has cyclic generator $1$ since $\lambda(R)\cdot 1 = R$.
The following provides necessary conditions for $A \| B$ to be a \textit{split
extension}, i.e.\ having $B \to A$ be a split monomorphism of the natural
$B$-$B$-bimodules, or a \textit{centrally projective extension},
i.e.\ the natural bimodule ${}_BA_B$ is centrally projective.  

\begin{lemma}
\label{lemma-sr}
If $A \| B$ is split, then ${}_SR$ is f.g.\ projective.
If $A \| B$ is centrally projective, then ${}_SR$ is a generator.
If $A \| B$ is split and centrally projective, then ${}_SR_{\mathcal{Z}(B)}$
is a faithfully balanced progenerator bimodule with $S$ and $\mathcal{Z}(B)$
Morita equivalent.  
\end{lemma}
\begin{proof}
In case $A \| B$ is split, apply the functor ${}_S\Hom(-, {}_BA_B)$ to ${}_BB_B \oplus * \cong {}_BA_B$,
obtaining ${}_SR \oplus * \cong {}_SS$. In case $A \| B$ is centrally
projective, apply the same functor to ${}_BA_B \oplus * \cong \oplus^n {}_BB_B$,
obtaining ${}_SS \oplus * \cong \oplus^n {}_SR$.  The last statement
follows, since $S \cong \End R_{\mathcal{Z}(B)}$ under the hypothesis
of centrally projectivity, and $\End {}_SR \cong \mathcal{Z}(B)$
follows from Proposition~\ref{lemma-below}, since centrally projective extensions
are depth two (or ``depth one implies depth two'') and
split extensions are right balanced (since $A_B$ is a generator,
${}_{\mathcal{E}}A_B$ is faithfully balanced).   
\end{proof}

%%%%%%%%%%%%%%%%%%%%%%%%%%%%%%%%%%%%%%%%%%%%%%%%%%%%%%%%%%%%%%%%%%%%%%%%%
\section{H-Separable Extensions}

H-separable extensions are useful to us as well-explored examples of
depth two extensions, ``toy models'' of depth two (D2) in a figure of speech.
Motivated by the question in \cite{KK} of whether D2 Hopf subalgebras
are normal, we first look at whether H-separable Hopf subalgebras are
normal Hopf subalgebras.  Recall that Hopf subalgebra $K$ of a Hopf
algebra $(H, m, u, \cop, \eps, \tau)$ is a Hopf algebra in its own right
w.r.t.\ the Hopf algebra structure (including the antipode $\tau$)
and is \textit{normal} if $\tau(a\1)K a\2 \subseteq K$ and
$a\1 K \tau(a\2) \subseteq K$ for all $a \in H$ (where $\cop(a) = a\1 \o a\2$).
Recall that a ring extension $A \| B$ is H-separable if the tensor-square
$A \o_B A$ is a centrally projective $A$-$A$-bimodule.  Equivalently,
in D2-friendly terms, this comes out as there being matched elements $e_i \in (A \o_B A)^A$
and $r_i \in R$  s.t.\ $\forall a,a' \in A$:
\begin{equation}
\label{eq: H-sep system}
a \o a' = \sum_i e_i \rho(r_i)(a) a' = \sum_i a \lambda(r_i)(a') e_i
\end{equation}
Note that $e_i \in T$ and $\lambda(r_i), \rho(r_i) \in S$.

\begin{lemma}[Sugano]
Let $A \| B$ be a $B$-projective H-separable extension.  Then
a two-sided ideal $I \lact A$ satisfies in terms of inducing and
contracting ideals:
$$ I = A( I \cap B) = (I \cap B)A.$$
\end{lemma}
\begin{proof}
Since $A_B$ is  projective, there is a dual basis $x_j \in A$
and $f_j \in \Hom(A_B, B_B) $, s.t.\ $\sum_j x_j f_j(a) = a$ $\forall a \in A$.
Since $A \| B$ is H-separable, we have $\End A_B \cong A \o_{\mathcal{Z}(A)} R^{\rm op}$
via left and right multiplication (similar to the Azumaya condition).  If $x \in I$,
then $x = \sum_j x_j f_j(x) \in A(I \cap B)$, since $f_j$ is $B$-valued right and left
multiplication by elements of $A$ on an ideal. Whence $I = A( I \cap B)$, and similarly
$I = (I \cap B) A$ by using dual bases for ${}_BA$. 
\end{proof}  

As its name indicates, an H-separable extension $A \| B$ is a \textit{separable extension}, i.e.\
the epi $\mu: A \o_B A \to A$ is split as an $A$-$A$-bimodule morphism.  We also see from the Azumaya
isomorphism mentioned in the proof, that $S \cong R \o_{\mathcal{Z}(A)} R^{\rm op}$.  It then follows
from Lemma~\ref{lemma-sr} that a split H-separable extension has a separable centralizer algebra over
$\mathcal{Z}(A)$, since ${}_SR = {}_{R^e}R$ is projective.  

\begin{theorem}
If $B \subseteq A$ is an H-separable  finite-dimensional Hopf subalgebra pair, then $B$
is a normal Hopf subalgebra of $A$.
\end{theorem}
\begin{proof}
It will suffice to show that the ideal $B^+ = \ker \eps|_B$ satisfies
$AB^+ = B^+ A$ by \cite[3.4.4]{Mo}.  First note that $B^+ = A^+ \cap B$,
the contracted ideal of the similarly defined $A^+ = \ker \eps$. 
We have the Nicholls-Zoeller theorem informing us that $A_B$ and ${}_BA$
are free f.g.\ modules.  Then by lemma $A^+ = AB^+ = B^+ A$.
\end{proof}

This method of showing normality extends only partially to the D2 case, where one is
limited to applying the method to $S$-stable ideals.  Although useful as a method of proof, the theorem
above is limited in applicability in some cases to trivial subalgebras $B = A$ as we see next.

\begin{prop}
Let $H < G$ be a subgroup of a finite group, and $A = \C G$, $B = \C H$ be the corresponding group algebras.
If $A \| B$ is H-separable, then $A = B$.  
\end{prop}
\begin{proof}
Let $\psi, \phi$ be irreducible characters on $G$.  The H-separability condition ${}_AA\o_B A_A \oplus * \cong
\oplus^n {}_AA_A$ becomes $\exists n \in \N:$ 
\begin{equation}
\label{eq: H-sep-char}
\langle {\rm Ind}_H^G {\rm Res}_H^G \psi \| \phi \rangle \leq n \langle \psi \| \phi  \rangle
\end{equation}
by recalling ${\rm Ind}\, V = A \o_B V$ for each simple module ${}_BV$.  
Let $\psi = 1_G$, the trivial character.  Then ${\rm Res}_H^G 1_G = 1_H$. It suffices
to show that $ {\rm Ind}_H^G 1_H = 1_G$ forcing $[G:H] = 1$ by dimensionality.  If $\phi \neq 1_G$,
then $\langle 1_G \| \phi \rangle = 0$ on the right, so the left-hand side becomes $\langle  {\rm Ind}_H^G 1_H \| \phi \rangle = 0$.
If $\phi = 1_G$, then using Frobenius reciprocity
$$\langle  {\rm Ind}_H^G 1_H \| 1_G \rangle = \langle 1_H \| {\rm Res}_H^G 1_G \rangle = 1.$$
Hence, $ {\rm Ind}_H^G 1_H = 1_G$.
\end{proof}
By switching from irreducible characters to simple modules and from inner products to hom-groups, we
may extend this triviality result to semisimple Hopf algebras.  

%%%%%%%%%%%%%%%%%%%%%%%%%%%%%%%%%%%%%%%%%%%%%%%%%%%%%%%%%%%%%%%%%%%%%%%%%%%%%%%%%%%%%%%%%%%%%%%%%%%%%%%%%%%%%%%%%%%%%
\section{Depth Two}

Recall that a \textit{depth two} ring extension $A \| B$ is characterized by its tensor-square
$A \o_B A$ being centrally projective w.r.t.\ the natural $B$-$A$-bimodule $A$ (left D2)
and the natural $A$-$B$-bimodule $A$ (right D2). Centrally projective ring extensions, H-separable
extensions and f.g.\
Hopf-Galois extensions are some of the classes of examples of D2 extension. If $A$ and $B$ are the complex
group algebras corresponding to a subgroup $H < G$ of  a finite group, then $A \| B$ is D2 iff $H$ is
a normal subgroup in $G$ \cite{KK}. 
 A new example to this list is the following:
\begin{prop}
\label{prop-wha}
Let $H$ be a finite-dimensional weak Hopf algebra, $A$ an $H$-comodule algebra, where coaction
$\rho_A(a) = a\0 \o a\1 \in A \o H$, and let $B = A^{\rm co\, H}$ be the subalgebra of coinvariants
where $\rho_A(b) = b1\0 \o 1\1$ for all $b \in B$.  If $A \| B$ is a weak $H$-Galois extension \cite{CDG}, then
$A \| B$ is right D2.
\end{prop}
\begin{proof}
Recall that $x := \rho_A(1) \neq 1_A \o 1_H$ necessarily, but that $x^2 = x$ and $x$ is a group-like
element in the Galois $A$-coring $\mathcal{C} := \Im g$ where $g: A \o H \to A \o H$ is a projection defined by
$g(a \o h) = a1\0 \o h 1\1$.  The $A$-coring structure on $\mathcal{C}$ is given by 
$a \cdot (a' 1\0 \o h 1\1) = aa' 1\0 \o h 1\1$, $ (a 1\0 \o h 1\1) \cdot a' = a {a'}\0 \o h {a'}\1$, coproduct 
$\cop_{\mathcal{C}}(a 1\0 \o h 1\1) = (a 1\0 \o h\1 1\1)\o_A (1 \o h\2)$ and
counit $\eps_{\mathcal{C}}(a 1\0 \o h 1\1) = a 1\0 \eps (h 1\1)$.    
The Galois structure is given by the isomorphism of $A$-corings ${\rm can}: A \o_B A \to \mathcal{C}$,
${\rm can}(a \o a') = axa' = a{a'}\0 \o {a'}\1$.  We note that ${\rm can}$ is an $A$-$B$-bimodule isomorphism, since 
 $b 1\0 \o 1\1 = 1\0 b \o 1\1$ for all $b \in B$ follows from \cite[(13),(16)]{CDG}. Since also $g$ is an $A$-$B$-bimodule
projection, $A \o_B A$ is isomorphic to a direct summand of $A \o H \cong \oplus^n A$ as $A$-$B$-bimodules, where
$n = \dim H$.  
\end{proof}
We note that $A \| B$ is a \textit{Frobenius extension} if
$A \o_B A \cong \mathcal{E}$ via a cyclic $A$-$A$-bimodule generator
$E: A \to B$.  For Frobenius extensions, right D2 is equivalent to left D2 \cite[6.4]{KS}.  

 In \cite{KS} we showed that $S := \End {}_BA_B$ is a left bialgebroid over $R$, and $T$ is a right bialgebroid over
$R$:  $S$ and $T$ are dual bialgebroids w.r.t.\ the nondegenerate pairing $\bra \alpha \| t \ket := \alpha(t^1)t^2$
with values in $R$ \cite{KS}. We also recall that $S$ acts on $A$ by evaluation with the invariant subring $A^S = B$
if $A_B$ is balanced.  
We continue next our study of the module ${}_SR$.  

\begin{prop}
\label{lemma-below}
Let $A \| B$ is left D2, then $\mathcal{E} \o_S R \cong A$ via evaluation.  If in addition $A_B$ is balanced,
then $\End {}_SR \cong \mathcal{Z}(B)$.  If $A \| B$ is a Frobenius extension, ${}_SR$ is a generator and
$\mathcal{E} \o_S R \cong A$ as $A$-$B$-bimodules, then $A$ is right D2.
\end{prop}
\begin{proof}
We have $A \o_R S \cong \mathcal{E}$ via $a \o \alpha \mapsto \lambda(a) \circ \alpha$ by \cite[3.10]{KS}.
Then $$ \mathcal{E} \o_S R \cong A \o_R S \o_S R \cong A \o_R R \cong A,$$
so that an inverse of $f \o r \mapsto f(r)$ is given by $a \mapsto \lambda(a) \o 1$. 

Next, if $f \in \End {}_SR$, then $\alpha \cdot f(r) = f(\alpha(r))$ for each $\alpha \in S$. 
Letting $\alpha = \lambda(r)$ where $r \in R$, we see that $f$ is determined by its value on $1$,
say $f(1) := a$, so $f = \rho(a)$.  Now $$ \alpha \cdot a = \alpha \cdot f(1) = f(\alpha(1)) = \alpha(1) a $$
since $\alpha(1) \in R$.  Whence $a \in A^S = B$, so $a \in B \cap R = \mathcal{Z}(B)$ and $f \in \lambda(\mathcal{Z}(B))$.

To prove the last statement, we note that ${}_SS \oplus * \cong \oplus^n {}_SR$ if ${}_SR$ is a generator.
Tensoring by ${}_A\mathcal{E}_B \o_S -$ and applying the hypothesis, we arrive at 
$${}_A \mathcal{E}_B \oplus * \cong \oplus^n {}_AA_B$$
If $A \| B$ is a Frobenius extension, then $\mathcal{E} \cong A \o_B A$ as $A$-$A$-bimodules, from which
we obtain that the tensor-square is centrally projective w.r.t.\ ${}_AA_B$.  
\end{proof}

We recall that equivalent conditions for depth two are the existence of (a left D2 quasibase) $\beta_i \in S$,
$t_i \in T$  and (right D2 quasibase) $\gamma_j \in S$, $u_j \in T$ such that
\begin{equation}
\label{eq: d2 quasibase}
a \o a' = \sum_i t_i \beta_i(a)a' = \sum_j a \gamma_j(a')u_j 
\end{equation}
for all $a, a' \in A$ (so by Eq.~(\ref{eq: H-sep system}) H-separability is D2).
We show that a split D2 extension is automatically f.g.\ projective.

\begin{prop}
\label{prop-fgp}
If $A \| B$ is a split D2 extension, then $A_B$ and ${}_BA$ are f.g.\ projective modules.  
\end{prop}
\begin{proof}
Suppose $p : A \to B$ is a conditional expectation, i.e.\ a bimodule projection onto $B$.  Applying
$p$ to the right tensorands of $a \o 1 = \sum_{i=1}^n t^1_i \o t^2_i \beta_i(a)$ yields
$$a = \sum_{i=1}^n t^1_i p(t^2_i \beta_i(a))$$ since $p(1) = 1$: a dual basis equation for $A_B$.
Whence $A_B$ is f.g.\ projective, and by a similar argument using a right D2 quasibase ${}_BA$ is
f.g.\ projective.  
\end{proof}

 The proposition allows a weakening of the hypotheses in \cite[Theorem 3.14]{KS}
that a biseparable D2 extension is automatically a \textit{QF-extension} $A \| B$, 
i.e.\ $A_B$ and ${}_BA$ are f.g.\ projective
and the left and right $B$-duals of $A$ are centrally projective w.r.t. the $A$-$B$-bimodule and $B$-$A$-bimodule $A$,
respectively.

%%%%%%%%%%%%%%%%%%%%%%%%%%%%%%%%%%%%%%%%%%%%%%%%%%%%%%%%%%%%%%%%%%%%%%%%%%%%%%%%%%%%%%%%%%%%%%%%%%%%%%%%%%%%%%%%
\section{A Depth Two Characterization of Semisimple-Hopf-Galois Extensions}

In this section we let $K$ denote a field of characteristic zero.  In this case, a finite dimensional Hopf algebra $H$
is semisimple if and only if it is cosemisimple by a theorem of Larson and Radford. 
 We limit ourselves here to semisimple-Hopf-Galois extension
$A \| B$, which means $H$ is a semisimple Hopf algebra, $A$ is a right $H$-comodule algebra, $B = A^{\rm co \, H}$,
and ${\rm can}: A \o_B A \to A \o H$ is an isomorphism. These are close to being classical Galois extensions
of noncommutative rings where groups of automorphisms are used in the definition.   

Thinking of the Steinitz characterization of finite degree field extensions that are Galois as being separable
and normal extensions, we see depth two filling the role of normal extension in the theorem below,
which is a clarification of parts of \cite[6.6]{KN} and \cite[8.14, 8.15]{KS}. An example of Frobenius extensions
with trivial centralizer are irreducible subfactors of finite Jones index such as the von Neumann algebras
of infinite-conjugacy-class group pairs.

\begin{theorem}
\label{th-char}
Suppose $A \| B$ is a Frobenius extension of $K$-algebras with trivial centralizer $R = 1_AK$.  
Then $A \| B$ is semisimple-Hopf-Galois if and only if $A \| B$ is a  separable and depth two extension.  
\end{theorem}
\begin{proof}
($\Rightarrow$) We have shown in Prop.~\ref{prop-wha} that $A \| B$ is right D2.  Similarly, $A \| B$ is left
D2 using an equivalent Galois mapping
${\rm can}': A \o_B A \to A \o H$ defined by ${\rm can}'(a \o a') = a\0 a' \o a\1$, also an isomorphism
since $\eta \circ {\rm can} = {\rm can}'$ where $\eta: A \o H \stackrel{\cong}{\longrightarrow} A \o H$,
$\eta(a \o h) = a\0 \o a\1 \tau(h)$. Finally an $H$-Galois extension where $H$ is semisimple is
 a separable Frobenius extension (Doi, Kreimer-Takeuchi).    

($\Leftarrow$) Since $A \| B$ is D2, we have dual bialgebroids $S$ and $T$ over $R$ \cite{KS}.  But $R$ is trivial, so
$S$ and $T$ are dual bialgebras over $K$.  Moreover, $S$ acts as a bialgebra
on $A$, $A_B$ is f.g.\ projective by the Frobenius hypothesis, and the
 right endomorphism ring is a smash product: $\mathcal{E} \cong A \# S$ via the mapping $a \# \alpha \mapsto
\lambda(a) \circ (\alpha \lact - )$. It will follow from \cite[8.3.3]{Mo} that $A \| B$ is a right Hopf-Galois extension
if we show $B = A^S$:
although the antipode plays no role, Schauenburg shows that $S$ and $T$ necessarily have antipodes, which
we can also see by showing they have nondegenerate integrals derived from the Frobenius structure $E: A \to B$
with dual bases $x_i, y_i \in A$.   

Since $A \| B$ is a separable Frobenius extension, there is $k \in K^*$ such that $\sum_i x_iy_i k = 1$.
But then $t_0 = \sum_i x_i \o y_i \in T$ is a nonzero integral with $\eps_T(t_0) = k^{-1}$,
 since $t_0 t = (t^1 t^2)t_0 = \eps_T(t) t_0$.  Whence $T$ is a semisimple Hopf algebra, therefore also
its dual $S$ is semisimple.  Consider $E \in S$, a nonzero integral since $\forall \alpha \in S$,
$\alpha E = \alpha(1) E = \eps_S(\alpha) E$.  Then $\eps_S(E) = E(1) = k'1_A \in K^* 1_A$, whence $A \| B$ is
split (with projection ${k'}^{-1}E$).  Then $A_B \cong B_B \oplus *$ is a $B$-generator, therefore balanced.
It follows from \cite[4.1]{KS} that the invariants in $A$ under the action of $S$ is $A^S = B$.  
We conclude that $A \| B$ is a $T$-Galois extension, where $T$ is a semisimple Hopf algebra.  
\end{proof}

%%%%%%%%%%%%%%%%%%%%%%%%%%%%%%%%%%%%%%%%%%%%%%%%%%%%%%%%%%%%%%%%%%%%%%%%%%%%%%%%%%%%%%%%%%%%%%%%%%%%%%%%%%%%%%%%%%%%%%%
\section{A Galois Coring}

We next show that a right balanced depth two extension $A \| B$ is a Galois
extension in the sense of Galois coring \cite{BW}. We continue
our notation for the right bialgebroid $T = (A \o_B A)^B$ over the
centalizer $R$, as well as right D2 quasibase $\gamma_j \in S$, $u_j \in T$
satisfying $\sum_j a\gamma_j(a')u_j = a \o_B a'$ for all $a,a' \in A$.

\begin{theorem}
\label{th-Galois coring}
If $A \| B$ is a D2 extension, then  $\mathcal{C} := A \o_R T$
with canonical structure is a Galois $A$-coring.  If $A_B$ is moreover
balanced or faithfully flat, then $B = A^{\rm co \, \mathcal{C}}$. 
\end{theorem}
\begin{proof}
First, there is a \textit{right $T$-comodule algebroid} structure $\rho: A
\to A \o_R T$ on $A$ given by
\begin{equation}
\label{eq: coaction}
\rho(a) = a\0 \o a\1 = \sum_j \gamma_j(a) \o u_j,
\end{equation}
with axioms we check next.  
That $\rho$ is coassociative and $\rho \o_R \id_T$ is
well-defined follows from using the isomorphism
$\beta: A \o_R T \stackrel{\cong}{\longmapsto} A \o_B A$
given by $\beta(a \o t) = at = at^1 \o t^2$ \cite[3.12(iii)]{KS}.  
Then $A \o_R T \o_R T \cong A \o_B A \o_B A$ via $\Phi := (\id_A \o \beta)(\beta \o \id_T)$,
so
$$\Phi( \id \o \cop_T) \circ \rho) = \sum_{j,k} \gamma_j(a) u^1_j \o_B \gamma_k(u^2_j)u^1_k \o_B u^2_k =
\sum_k 1 \o \gamma_k(a)u^1_k \o u^2_k = 1 \o 1 \o a
$$
$$= \sum_{j,k} \gamma_k(\gamma_j(a))u^1_k \o u^2_ku^1_j \o u^2_k = 
\Phi( (\rho \o \id)\rho(a)). $$ 
Also, $a\0 \eps_T(a\1) = $ $ \sum_j \gamma_j(a)u^1_j u^2_j  = a$
for all $a \in A$.  

For $\rho(aa') = \rho(a)\rho(a')$ to make sense in an appropriate tensor subalgebra in
$A \o_R T$ we need to check:
$$ \beta (r \cdot a\0 \o a\1) = \sum_j r\gamma_j(a)  u_j = \sum_j
\sum_j \gamma_j(a)u^1_j r \o u^2_j = \beta(a\0 \o \tilde{t}(r) a\1). $$
where antihomomorphism $\tilde{t}(r) = 1 \o r$, $ R \to T$ is the target map. Next,
$$\beta(\rho(a)\rho(a')) = \sum_{j,k} \gamma_j(a) \gamma_k(a') u_j u_k =
 \sum_{j,k} \gamma_j(a) \gamma_k(a') u_k^1 u_j^1 \o u_j^2 u_k^2 
$$
$$ = 1 \o aa' = \sum_j \gamma_j(aa')u_j = \beta(\rho(aa')). $$
Also $\rho(1_A) = 1_A \o_R 1_T$ since $\gamma_j(1_A) \in R$.
Finally we note that for each $b\in B$
$$ \rho(b) = \sum_j \gamma_j(b) \o_R u_j = b \o \sum_j \gamma_j(1)u_j = b \o 1_T $$
so $B \subseteq A^{\rm co \, \rho}$.  
The converse: if $\rho(x) = x \o 1_T$ $ = \sum_j \gamma_j(x) \o u_j$ applying $\beta$ we
obtain $x \o_B 1 = 1 \o_B x$.  If $A_B$ or ${}_BA$ is  faithfully flat, $x \in B$.
If $A_B$ is balanced, we know $A^S = B$ under the action $\lact$ of $S$ on $A$ \cite[4.1]{KS}.
Applying $\mu (\alpha \o \id)$ for
each $\alpha \in S$, we obtain $\alpha \lact x =$ $ \alpha(x) =$ $ \alpha(1) x$,
whence $x \in B$ also in this case.
 
A left and right $A$-action $\MC$ is then given by
\begin{equation}
\label{eq: A-action}
 a \cdot (a' \o t) = aa' \o t \ \ \ (a \o t) \cdot a' = a{a'}\0 \o t {a'}\1 
\end{equation}
which satisfies $$(a \o t) \cdot 1_A = \sum_j a \gamma_j(1) \o_R tu_j =
a \o \sum_j \gamma_j(1)u^1_j t^1 \o t^2 u_j^2 = a \o t .$$
The rest of the axioms of an $A$-$A$-bimodule structure on $\MC$ follow readily.
A comultiplication on $\MC$ is given by
\begin{equation}
\label{eq: coring coprod}
\cop_{\MC}(a \o t) = (a \o t\1) \o_A (1_A \o t\2)
\end{equation}
in terms of $\cop_T(t) = \sum_j t^1 \o_B \gamma_j(t^2) \o_R u_j$ \cite[(82)]{KS}
i.e., $\cop_{\MC}(a \o t) = a \o t\1 \o t\2$ $\in A \o_R T \o_R T \cong \MC \o_A \MC$.  
This coproduct is clearly left $A$-linear, and right $A$-linear since $\rho$ is coassociative: 
$$\cop_{\MC}(a \o t) \cdot a' = (a \o t\1) \o_A ({a'}\0 \o t\2 {a'}\1) =
(a{a'}\0 \o t\1 {a'}\1 ) \o (1_A \o t\2 {a'}\2) = \cop_{\MC}((a \o t) \cdot a') .
$$
$\cop_{\MC}$ is coassociative since $\cop_T$ is so. 

The counit of $\MC$ is given by 
\begin{equation}
\label{eq: coring counit}
\eps_{\MC}(a \o t) = a \eps_T(t) = a t^1 t^2,
\end{equation}
clearly left $A$-linear.  The counit is also right $A$-linear since
$$ \eps_{\MC}((a \o t)\cdot a') = a {a'}\0 \eps_T(t {a'}\1) = \sum_j a \gamma_j(a') u_j^1t^1 t^2 u_j^2  
 = a t^1 t^2 a' = \eps_{\MC}(a \o t) a'.$$
The counitality axioms follow from
$$ \eps_{\MC}(a \o t\1) 1_A \o_R t\2 = \sum_j a t^1 \gamma_j(t^2) \o_R u_j = a \o t ,$$
$$ a \o t\1 \eps_{\MC}(1 \o t\2) = a \o_R (t^1 \o_B \sum_j \gamma_j(t^2)u^1_j u^2_j) = a \o t.$$

The Galois structure on $\MC$ comes from the grouplike element $x := 1_A \o 1_T$
($\cop_{\MC}(x) = x \o_A x$ since $\cop_T(1_T) = 1_T \o_R 1_T$,
and $\eps_{\MC}(x) = 1$ since $\eps_T(1_T) = 1$).  The Galois canonical mapping ${\rm can}: A \o_B A \to \MC$
given by 
\begin{equation}
\label{eq: can}
{\rm can}(a \o a') = a x a' = a {a'}\0 \o_R {a'}\1
\end{equation}
 is well-defined since $B \subseteq A^{\rm co \, \rho}$ or
$x \cdot b = b \cdot x$ for each $b \in B$. It is an isomorphism (of $A$-corings where $A \o_B A$
is armed with its Sweedler canonical structure with grouplike element $1 \o 1$) since ${\rm can} \circ\beta = \id$
and $\beta \circ {\rm can} = \id$.  
\end{proof}

We propose to  consider the following set-up in a future paper.  
Suppose $T$ is a right $R$-bialgebroid where ${}_RT$ and $T_R$ are f.g.\ projective
with unique dual left $R$-bialgebroid $S$.
Define a (right) \textit{$T$-Galois extension} $A \| B$ to be a right $T$-comodule algebroid $A$
such that $A \o_R T$ is a Galois $A$-coring with structure given by Eqs.~(\ref{eq: A-action}),~(\ref{eq: coring coprod}),~(\ref{eq: coring counit})
and~(\ref{eq: can}) and $B = A^{\rm co \, T}$.   
Then supposing $A \| B$ to be a Frobenius extension of rings, we will show 
$A \| B$ is a $T$-Galois extension if and only if $A \| B$ is depth 2 and balanced.  
The theorem has established the implication $\Leftarrow$: we note that $T_R$ and ${}_RT$ are f.g.\ projective from specializing 
the left and right quasibases equations to $T$. The implication ``Galois $\Rightarrow$ D2'' follows from \cite[2.1(4)]{KK} transposed
to: $A \o_R T \cong A \o_B A$ as $A$-$B$-bimodules and ${}_RT$ f.g.\ projective iff $A \| B$ is right D2. Since
$A \| B$ is Frobenius, it is then also left D2.

%%%%%%%%%%%%%%%%%%%%%%%%%%%%%%%%%%%%%%%%%%%%%%%%%%%%%%%%%%%%%%%%%%%%%%

\end{document}